\documentclass[12pt, a4paper]{article}
\usepackage{cmap}      % search in PDF
\usepackage{mathtext}  % russian letters in formulas
\usepackage[T2A]{fontenc}   % encoding
\usepackage[utf8]{inputenc} % source text encoding
\usepackage[russian,english]{babel} % localization and hyphenation (English set as primary)
\usepackage{amsfonts}       % mathematical symbols
\usepackage{graphicx}       % pictures
\usepackage{amsmath}
\usepackage{MnSymbol}
\usepackage{wasysym}
\graphicspath{{pictures/}}  % pictures directory
\DeclareGraphicsExtensions{.jpg}    % formats
\title{\textbf{Mathematical Billiards}}
\author{Grigoriy Yakovlev}
\date{}

\begin{document}

\maketitle

\begin{abstract}

This paper considers planar mathematical billiards with smooth billiard tables. We prove that through any point outside a convex smooth caustic, there passes a unique billiard table (Lemma 2). Furthermore, we describe a billiard table having a nephroid as its caustic.

\end{abstract}

\section{Introduction}

\quad

Everyone is familiar with the game of billiards on a rectangular table with pockets. There are many variations of this game. Usually, they differ from each other in rules and the number of balls. For example, French billiards has no pockets, which makes it similar to mathematical billiards. A mathematical billiard consists of a certain region (the billiard table) and a point that moves freely---that is, rectilinearly and at a constant speed---within this region. Upon reaching the boundary, the point reflects off it according to the law of geometric optics: the angle of incidence equals the angle of reflection. After reflection, the rectilinear motion continues until the next collision with the boundary (see Fig. 1).

\begin{figure}[h]
    \centering
    \includegraphics[width=70mm]{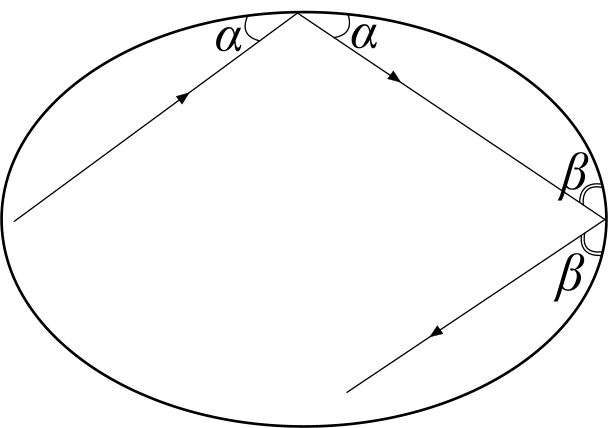}
    \caption{\textit{Mathematical billiards is a dynamical system}}
\end{figure}

Consider a billiard in an ellipse. An ellipse is a figure on the plane defined by the equation $\frac{x^2}{a^2}+\frac{y^2}{b^2}=1$. One can consider an ellipse as a cross-section of a cylinder or a circle stretched along one of its axes. The ellipse itself possesses interesting properties. We are currently interested in the focal property of the ellipse: an ellipse is the geometric locus of points such that the sum of the distances from them to two given points, called \textit{foci}, is a constant. Ellipses that share the same foci are called \textit{confocal}.
\quad

\textbf{Theorem (The billiard property of the ellipse)}. If a point moves along a tangent to a certain confocal ellipse, then after reflection it will also move along another tangent to the same confocal ellipse (see Fig. 2).

That is, if one segment of the billiard trajectory is tangent to a confocal ellipse, all subsequent segments will also be tangent to it. The proof of this theorem can be found in \cite{GAL}.

\begin{figure}[h]
    \centering
    \includegraphics[width=70mm]{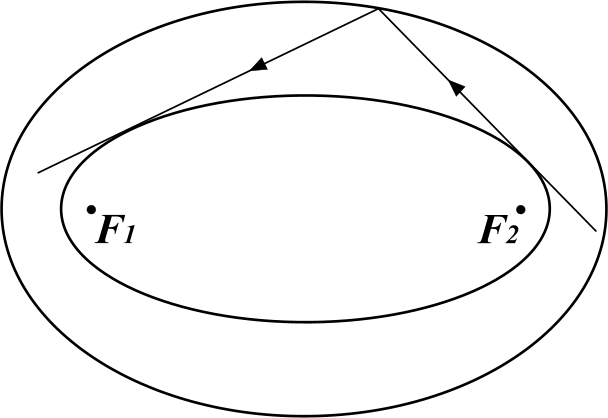}
    \caption{\textit{The billiard property of the ellipse. The figure shows two confocal ellipses and two segments of a billiard trajectory that are tangent to the smaller ellipse.}}
\end{figure}

This remarkable observation leads us to the following definition. Let $\Gamma$ be a billiard table, and $\gamma$ be a smooth curve inside this billiard table such that if our trajectory touches it at least once, then upon every subsequent reflection, the trajectory will also be tangent to it. In this case, $\gamma$ is called a \textit{caustic} for our billiard table $\Gamma$.

It turns out that, other than the ellipse, there are no known examples of billiard tables with a caustic that are algebraic curves. Our work consisted of searching for such a billiard table-caustic pair.

\quad

There is a method that allows constructing a billiard table from a given caustic: if we take a smooth convex closed curve $\gamma$, stretch an inextensible string of a fixed length $L$ around it, place a pencil at point $B$ (see Fig. 3), and make a revolution around our curve, we will draw a billiard table for which $\gamma$ will be a caustic. This method is called the \textit{string construction}. Its justification can be found in \cite{TAB}.

\begin{figure}[h]
    \centering
    \includegraphics[width=70mm]{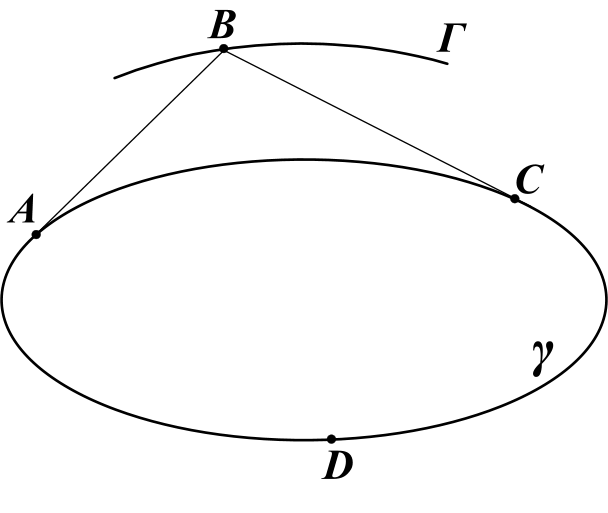}
    \caption{\textit{String construction}}
\end{figure}

In practice, it is quite difficult to use this method to find a billiard table from a caustic. The difficulty lies in the fact that, since the string is inextensible, one must constantly ensure that the value $AB+BC+\smile\!\!ADC$ remains constant at every moment in time. The arc length $\smile\!\!ADC$ is calculated using the following formula:

\begin{center}
$$\smile\!\!ADC=\int\limits_{t_1}^{t_2}\sqrt{\bigl(x'(t)\bigl)^2+\bigl(y'(t)\bigl)^2}dt.$$
\end{center}

Where $(x(t),y(t))$ is the parameterization of the curve $\gamma$, and $t_1$ and $t_2$ are the parameter values. This integral cannot always be expressed in terms of elementary functions. For example, if the caustic is an ellipse, the required arc length is expressed only using special functions.

\quad

But even if the arc length is expressed in elementary functions, this does not necessarily mean that we can explicitly write down the equation of the billiard table. The equation of the billiard table is given in the form $F(x,y,L)=0$, where $L$ is the parameter corresponding to the string length. For example, in the case of a circle, the function $F$ takes the form
$$F(x,y,L)=x^2+y^2-\sec^2(t),$$
where $t$ is the smallest positive root of the equation $\tan(t)-t=\frac{L}{2}-\pi$. This equation is unsolvable in elementary functions, as shown in \cite{KAN}.

From geometric considerations, we can establish that if we apply the string construction to a circle, we will obtain a billiard table in the form of a concentric circle of a larger radius.

\begin{figure}[h]
    \centering
    \includegraphics[width=90mm]{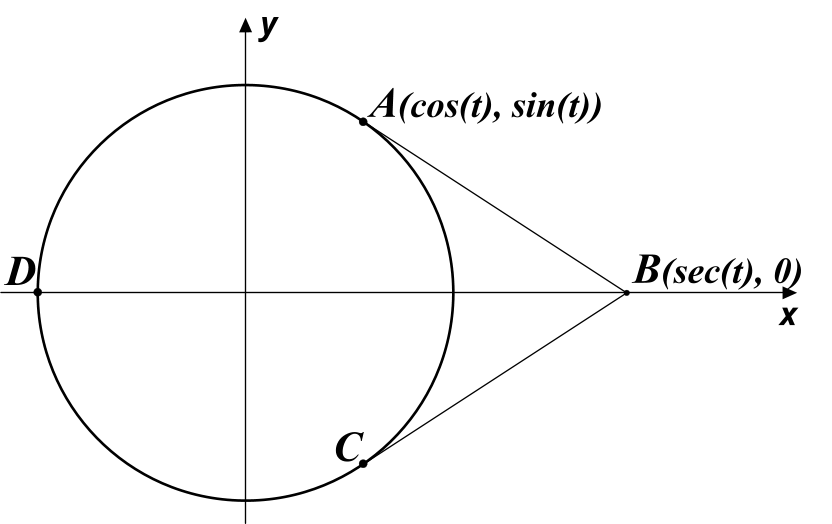}
    \caption{\textit{Circle}}
\end{figure}

\quad

Let us assume that our circle has a radius of 1, and the string length is $L$. We center this circle at the origin of the coordinate system. Points on this circle can be parameterized as follows:

\begin{center}
$\begin{cases}
x(t)=\cos(t);\\
y(t)=\sin(t).\\
\end{cases}$
\end{center}

Let $B$ be the point where "the string is stretched by the pencil" (see Fig. 4). By symmetry, we can assume that point $B$ lies on the $Ox$ axis. Then our task is reduced to finding the $x$-coordinate of point $B$, which will be the radius of the resulting circle. 

Let $A\bigl(\cos(t),  \sin(t)\bigl)$ be one of the points where the string detaches from the caustic. The general form of the tangent equation to a curve with parameterization $(x(t),y(t))$ at a point $t_0$ is $\bigl(x - x(t_0)\bigl)y'(t_0) - \bigl(y - y(t_0)\bigl)x'(t_0) = 0$. Substituting the point $A\bigl(\cos(t), \sin(t)\bigl)$ into this equation, we get: $x \, \cos(t) + y \, \sin(t) = 1.$ Hence, the $x$-coordinate of point $B$ is $\sec(t)$---this will be the radius of the larger circle. As we can see, to find the radius of the billiard table, we need to find $t$, which can be found from the following relation:
$\pi - t + \sqrt{(\cos(t)-\sec(t))^2+\sin^2(t)} = \frac{L}{2}$, since $\smile\!\!DA + AB = \frac{L}{2}$, where $\smile\!\!DA$ is the smaller of the two arcs. Simplifying this expression yields: $\tan(t) - t = \frac{L}{2} - \pi$.  But since we cannot express $t$ in terms of elementary functions, we also cannot express the equation of the billiard table, as it has the form $x^2 + y^2 = \sec^2(t)$.

\quad

Now let us consider the curve with the following parameterization:

\begin{center}
$\begin{cases}
x(t) = \sin(t) + \frac{1}{3} \sin(3t);\\
y(t) = \cos(t) + \frac{1}{3} \cos(3t).\\
\end{cases}$
\end{center}

This curve is called a nephroid, and it is a special case of an epicycloid (see Fig. 5). Although it is neither smooth nor convex, we can take its upper half and paste the segment $[-\frac{2}{3}; \frac{2}{3}]$ into the gap (see Fig. 6); in this case, the curve becomes convex, and we can apply the string construction to it. It should be noted that the curve will not be infinitely smooth, since at the junction points of the nephroid with the segment, the curve will only be once differentiable. Therefore, we consider only the arc of the billiard table obtained using the string construction during the time intervals when the string completely covers the rectilinear segment of our curve.

\begin{figure}[h!]
\centering
\begin{minipage}{.4\textwidth}
  \centering
  \includegraphics[scale=0.425]{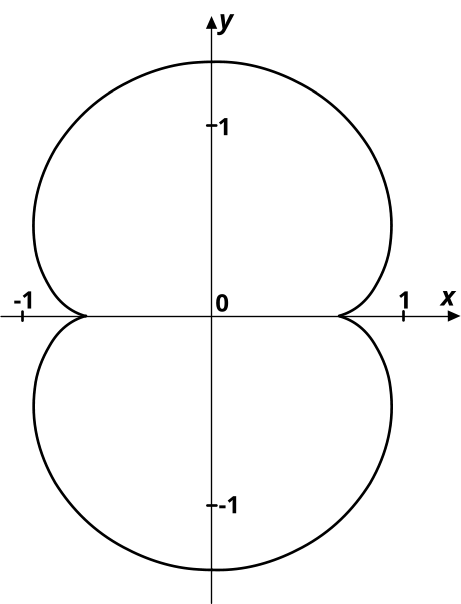}
  \caption{Nephroid}
\end{minipage}%
\begin{minipage}{.7\textwidth}
  \centering
  \includegraphics[scale=0.45]{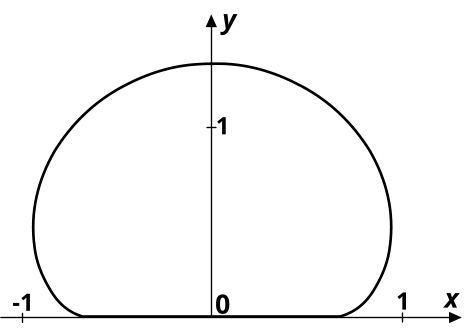}
  \caption{Upper half of the nephroid}
\end{minipage}
\end{figure}

Also, the nephroid is indeed an algebraic curve, as it is defined by the following polynomial of the sixth degree:

\begin{center}
$(x^2+y^2-\frac{4}{9})^3=\frac{4}{3}y^2.$
\end{center}

Therefore, if the arc of the billiard table we are interested in turns out to be algebraic, we will be able to find a billiard table-caustic pair that possesses corresponding algebraic arcs.

\quad

\textbf{Theorem.} The considered arc of any billiard table for which the upper half of a nephroid with an inserted segment serves as a caustic is an algebraic curve.

\section{Proof of the Theorem}

\quad

The arc length of the curve under consideration is calculated by the formula $2 \bigl(\sin(t_{2}) - \sin(t_{1})\bigl)$, where $t_{1}$ and $t_{2}$ are parameter values. As we can see, the arc length is expressed only in terms of trigonometric functions of $t$; therefore, when calculating the string length, we should not encounter the problems that arose with the circle, where we obtained both a trigonometric function of $t$ and $t$ itself in a single equation.

Then, substituting $t_{1}$ and $t_{2}$ into the tangent equation and solving this system of equations for $x$ and $y$, we obtain the intersection point of the tangents drawn to the curve at the points corresponding to parameters $t_{1}$ and $t_{2}$:

\begin{center}
$\begin{cases}
x = \frac{2}{3}\bigl(\cos(t_{1}-2t_{2})-\cos(2t_{1}-t_{2})-\cos(2t_{1}+t_{2})+\cos(t_{1}+2t_{2})\bigl)\csc\bigl(2(t_{1}-t_{2})\bigl);\\
y = \frac{8}{3}\bigl(\sin(t_{1})-\sin(t_{2})\bigl)\cos(t_{1})\cos(t_{2})\csc\bigl(2(t_{1}-t_{2})\bigl).\\
\end{cases}$
\end{center}

That is, $x$ and $y$ are rational functions of $z:=e^{i t_{1}}$ and $w:=e^{i t_{2}}$.

\quad

Now, knowing the arc length formula and the intersection point of the tangents, we can calculate the length of the string that detaches from our caustic at the considered points. It turns out that this is also a rational function of $z$ and $w$, but some of the expressions will be inside an absolute value. To avoid case analysis, we prove the following lemma.

\quad

\textbf{Lemma 1.} If a surface is defined parametrically as $\gamma=\left\{(x(t),y(t),z(t))\right\},t \in\mathbb{R}^2$, where $x,y,z:\mathbb{R}^2\longrightarrow\mathbb{R}$ are rational functions, then $\gamma$ is contained in an algebraic surface.

\quad

\textbf{Proof.} Consider a polynomial $F(x,y,z)$ of three variables of degree $N$. This polynomial has $K_1(N):=\bigl({C_{N+1}^{2}}+(N+1)\bigl)+\bigl({C_{N}^{2}}+N\bigl)+\dots=N^2+1+\frac{N(N^2+11)}{6}$ distinct monomials. We want to find such values for the coefficients of these monomials that our polynomial identically equals zero as a function of $t$; this will mean that the resulting surface is contained in some algebraic surface. Let $t := (t_1, t_2)^{T}$. Let us also denote:

\begin{center}
$\begin{cases}
x:=\frac{P_1}{Q_1};\\
y:=\frac{P_2}{Q_2};\\
z:=\frac{P_3}{Q_3}.\\
\end{cases}$
\end{center}

Where $P_i,Q_i$ are polynomials in $t_1,t_2$.

Consider $F$ as a function of $t_1, t_2$, having previously multiplied it by $(Q_1Q_2Q_3)^N$ so that $F$ becomes a polynomial. The degree of this polynomial in two variables is bounded by $N(m+2n)$, where $m:=\max(\deg(P_{i})),n:=\max(\deg(Q_{i}))$. From this, we can conclude that the total number of distinct monomials in $F(t_1,t_2)$ will be at most $K_2(N):=\frac{(N(m+2n)+1)(N(m+2n)+2)}{2}$.

Since we want the coefficients of the monomials in $F(t_1,t_2)$ to vanish, we need to solve a system with $K_1$ variables and $K_2$ equations. Note that for sufficiently large $N$, the inequality $K_1>K_2$ holds, and since the system will be linear and homogeneous, it will have a nontrivial solution.

\begin{flushright}
$\blacksquare$
\end{flushright}

Returning to our problem. As a consequence of the lemma, there is an algebraic relation between $x, y$, and $L$ for each case of resolving the absolute values in the function $L$. Then, if we fix the value of $L$, we will obtain a polynomial $F(x,y)$ of two variables. Thus, we have proven the algebraicity of each of the arcs of the billiard table resulting from different resolutions of the absolute values.

\quad

Let us show that, in fact, the entire arc of the billiard table we are interested in is a single algebraic curve. To prove this fact, we will prove two auxiliary lemmas.

\quad

\textbf{Lemma 2.} If the caustic of a billiard table is a closed convex smooth curve with non-zero curvature, then the billiard table is also a smooth curve. Moreover, for any such caustic, through every point lying outside this curve, there passes exactly one billiard table for which this curve will be the caustic.

\quad

\textbf{Proof.} Let our curve be defined parametrically: $\gamma=\left\{(x(t),y(t))\right\}$,
$t \in\mathbb{R}$. Consider an arbitrary $\tau\in\mathbb{R}$ and such an open arc $\gamma_0$ of the given curve that $(x(\tau),y(\tau))\in\gamma_0$, and that the tangent vectors at any two distinct points of this arc are not collinear. This arc generates an open neighborhood $M \subset \mathbb{R}^2$, which corresponds to pairs of points from $\gamma_0$. We want to recover the intersection point of the tangents to our curve at the points $\gamma(t_1)$ and $\gamma(t_2)$. The coordinates $f_1, f_2$  of this point are found from the following system of equations with respect to $x, y$:

\begin{center}
$\begin{cases}
(x-x(t_1))y'(t_1)-(y-y(t_1))x'(t_1)=0;\\
(x-x(t_2))y'(t_2)-(y-y(t_2))x'(t_2)=0.
\end{cases}$
\end{center}

We need to check that this system has a unique solution. To do this, we calculate the determinant of the corresponding matrix: it will be equal to $x'(t_1)y'(t_2)-x'(t_2)y'(t_1)$. Without loss of generality, we can assume that our curve has an arc-length parameterization; then the velocity vector has a constant length everywhere and thus never vanishes. By rotating our curve and choosing a sufficiently small arc, we can ensure that $x'(t_0),y'(t_0)$ do not vanish for any point $t_0\in\gamma_0$. And then, since the tangent vectors to our curve at distinct points of the given arc cannot be collinear by construction, this expression cannot equal zero.

\quad

Note that we have obtained a smooth mapping $F:(t_1,t_2)\mapsto(f_1,f_2)$, since $x(t),y(t)$ are smooth functions of one variable. We want to construct a smooth mapping $F^{-1}:(f_1,f_2)\mapsto(t_1,t_2)$. To do this, we will use the Inverse Function Theorem. We must verify that the corresponding Jacobian does not vanish at a certain fixed $(a, b)\in M$. This is equivalent to the fulfillment of the following conditions:

\begin{center}
$\begin{cases}
x'(a) \big( y(a) - y(b) \big) \neq y'(a) \big(x(a) - x(b) \big);\\
x'(b) \big( y(a) - y(b) \big) \neq y'(b) \big(x(a) - x(b) \big);\\
x'(a)y''(a) \neq x''(a)y'(a);\\
x'(b)y''(b) \neq x''(b)y'(b).\\
\end{cases}$
\end{center}

Since the velocity vector can't have vanishing components on $\gamma_0$, $y(t_1)-y(t_2)$ and $x(t_1)-x(t_2)$ cannot be equal to zero on our arc. Then the first condition can be rewritten as $\frac{x'(a)}{y'(a)}\neq\frac{x(a) - x(b)}{y(a) - y(b)}$. Since the curve is  strictly convex, the vector $(x(a),y(a))-(x(b),y(b))$ cannot be parallel to $\gamma'(a)$ or $\gamma'(b)$, meaning that the first two conditions are indeed satisfied.

Since, by the condition of the lemma, the curvature at every point is non-zero, the arc $\gamma_0$ can be chosen such that $x''(t)$ and $y''(t)$ do not vanish on it. Then the third condition can be rewritten as $\frac{x'(a)}{y'(a)}\neq\frac{x''(a)}{y''(a)}$. Consider another curve defined by the following parameterization: $\gamma'=\left\{(x'(t),y'(t))\right\},t \in\mathbb{R}$. Since the curve $\gamma$ has an arc-length parameterization, $\gamma'$ is actually a circle, which means $(x'(t),y'(t))\perp(x''(t),y''(t))$ for any $t\in\mathbb{R}$. That is, these vectors cannot be collinear at any point of our curve, which means the third and fourth conditions are also satisfied.

Thus, there indeed exists a smooth mapping $F^{-1}:(f_1,f_2)\mapsto(t_1,t_2)$ acting from a certain open neighborhood of the point $F(a,b)$ and assigning to every point of this neighborhood the parameter values $t$, which, when substituted into the curve equation, give us the coordinates of the tangency points when drawing tangents to our curve from the chosen point. Then the vectors $v_1=(x(t_1)-f_1,y(t_1)-f_2);v_2=(x(t_2)-f_1,y(t_2)-f_2)$ will also smoothly depend on $(f_1,f_2)$.

\quad

Suppose a billiard table, defined in this neighborhood by the function $y=f(x)$ for which the given curve acts as a caustic, passes through an arbitrary point $(x_0,y_0)$ from the considered neighborhood. Then the segments of the tangents from this point to the caustic must form equal angles with the corresponding tangent to the table. This means that the direction vector of the tangent to the billiard table at the given point is uniquely determined up to multiplication by a non-zero constant. The normal of the considered tangent is expressed as $(h_1,h_2):=\frac{v_1}{|v_1|}+\frac{v_2}{|v_2|}$, and hence the direction vector has the form $(-h_2,h_1)$. From this we obtain the following condition for our table: $f'(x)=-\frac{h_1}{h_2}$.

Note that we have obtained an ordinary differential equation, the right-hand side of which is a smooth function of $x$ and $y$ on some open set. Hence, according to the existence and uniqueness theorem for solutions of ordinary differential equations, the proof of which can be found in \cite{PONT}, there exists a unique function $f$ that will satisfy this relation in a neighborhood $M$ of the point $(x_0,y_0)$, and the obtained function will be smooth. Thus, through any point in $M$, there passes a unique smooth arc of a billiard table for which $\gamma$ will be a caustic. Due to the arbitrariness in choosing the initial point $\tau\in\gamma$, exactly one billiard table for which $\gamma$ is a caustic will indeed pass through any point outside $\gamma$, and this table must be smooth.

\begin{flushright}
$\blacksquare$
\end{flushright}

\textbf{Corollary 1.} Using the string construction, we can obtain any billiard table for which a given curve satisfying the restrictions of the lemma will be a caustic.

\quad

\textbf{Corollary 2 (Graves' confocal ellipses).} If we take an ellipse and apply the string construction to it, we will again obtain an ellipse.

\quad

\textbf{Lemma 3.} Let an infinitely smooth curve $\gamma$ be piecewise defined by irreducible polynomials:
$F_{1}(x,y)=0$ for $x\leq x_{0}$,
$F_{2}(x,y)=0$ for $x\geq x_{0}$,
$F_{1}(x_{0}, y_{0})=F_{2}(x_{0}, y_{0})$.
Then $F_{1}(x, y)= \lambda F_{2}(x, y)$, $\lambda \in \mathbb{R}\setminus \left\{ 0 \right\}$.

\quad

\textbf{Proof.} Let us make a coordinate transformation such that the partial derivatives of $F_{1}, F_{2}$ with respect to $y$ do not vanish at the point $(x_{0},y_{0})$. Then, according to the Implicit Function Theorem for analytic functions, the proof of which can be found in \cite{GRIF}, there exists a certain neighborhood $M(x_{0}, y_{0})$ on which analytic functions $f_{1}(x), f_{2}(x)$ exist such that: 

\begin{center}
$\begin{cases}
F_{1}(x, y)=0 \Leftrightarrow y=f_{1}(x);\\
F_{2}(x, y)=0 \Leftrightarrow y=f_{2}(x).\\
\end{cases}$
\end{center}

Then $F_{1}$ and $F_{2}$ can be written on $M$ in the following form:
\begin{center}
$\begin{cases}
F_{1}(x, y)= \big( y - f_{1}(x) \big) G_{1}(x,y);\\
F_{2}(x, y)= \big( y - f_{2}(x) \big) G_{2}(x,y).\\
\end{cases}$

$G_{1}, G_{2} \neq 0$ $\forall (x,y) \in \mathbb{R}.$
\end{center}

\quad

Since $f_{1}, f_{2}$ are analytic, they can be expanded into a Taylor series in the neighborhood $M$. Since $\gamma$ is infinitely smooth, all derivatives of $f_{1}, f_{2}$ coincide at the point $(x_{0}, y_{0})$, meaning their Taylor series expansions also coincide, i.e., $f_{1}=f_{2}$. Then the polynomials $F_{1}$ and $F_{2}$ vanish on the same infinite set of points; consequently, the relation $F_{1}= \lambda F_{2}$ holds.

\begin{flushright}
$\blacksquare$
\end{flushright}

Returning to our problem. As a consequence of the second lemma, we can conclude that the considered arc of our billiard table, obtained for the nephroid using the string construction, will be smooth. Further, since each section of the table resulting from different resolutions of the absolute values is defined by a polynomial of two variables, and the table is smooth at the junction points, by the third lemma, the entire table is actually defined by a single algebraic equation. Also, it follows from the second lemma that using the string construction we can obtain any billiard table for the nephroid, which means that indeed, any billiard table for which the nephroid is a caustic will be defined by an algebraic equation.

\begin{flushright}
$\blacksquare$
\end{flushright}

{}

\end{document}